\documentclass[11pt,a4,reqno]{amsart}
\allowdisplaybreaks

\usepackage{amsfonts}
\usepackage{amsmath}
\usepackage{amssymb}
\usepackage{amsthm}
\usepackage{xcolor}
\usepackage{enumerate}
\usepackage{graphicx}
\usepackage{hyperref}
\usepackage{mathrsfs}
\usepackage{stmaryrd}
\usepackage[nameinlink,capitalize]{cleveref}

\newcommand{\be}{\begin{eqnarray}}
\newcommand{\ee}{\end{eqnarray}}
\newcommand{\ce}{\begin{eqnarray*}}
\newcommand{\de}{\end{eqnarray*}}

\newtheorem{theorem}{Theorem}[section]
\newtheorem{proposition}[theorem]{Proposition}
\newtheorem{lemma}[theorem]{Lemma}
\newtheorem{corollary}[theorem]{Corollary}
\newtheorem{remark}[theorem]{Remark}
\newtheorem{definition}[theorem]{Definition}
\newtheorem{example}[theorem]{Example}
\newtheorem{assumption}[theorem]{Assumption}

\def\bt{\begin{theorem}}
\def\et{\end{theorem}}
\def\bp{\begin{proposition}}
\def\ep{\end{proposition}}
\def\bl{\begin{lemma}}
\def\el{\end{lemma}}
\def\bc{\begin{corollary}}
\def\ec{\end{corollary}}
\def\bd{\begin{definition}}
\def\ed{\end{definition}}
\def\br{\begin{remark}}
\def\er{\end{remark}}
\def\bx{\begin{example}}
\def\ex{\end{example}}
\def\ba{\begin{assumption}}
\def\ea{\end{assumption}}

\def\[{{\Big[}}
\def\]{{\Big]}}
\def\<{{\langle}}
\def\>{{\rangle}}
\def\({{\Big(}}
\def\){{\Big)}}

\def\>{\geqslant}
\def\<{\leqslant}

\def\bD{{\mathbf D}}

\def\cF{{\mathcal F}}
\def\cG{{\mathcal G}}
\def\cH{{\mathcal H}}

\def\cN{{\mathcal N}}

\def\mB{{\mathbb B}}

\def\mD{{\mathbb D}}
\def\mE{{\mathbb E}}

\def\mP{{\mathbb P}}

\def\mR{{\mathbb R}}

\def\Om{\Omega}

\def\esssup{{\rm ess sup}}
\newcommand{\sptext}[3]{\hspace{#1 em}\mbox{#2}\hspace{#3 em}}

\begin{document}


\title{A note on the Besov spaces $\mB_1^{\Phi_2}$}

\author{Stefan Geiss}
\address{Department of Mathematics and Statistics, University of Jyv\"askyl\"a, Finland}
\email{stefan.geiss@jyu.fi}
\date{\today}
\begin{abstract}
Regarding the Besov space $\mB_1^{\Phi_2}$, defined by S. Geiss and Ylinen (Memoirs AMS 1335, 2021)
by a coupling method, we prove for $\xi\in\mB_1^{\Phi_2} \cap L_2$ a distributional estimate 
in terms of the Lorentz-Zygmund space $L_{1,\infty}(\ln L)^{-1}$ for the
Malliavin derivative $\bD(T_u\xi)$ uniformly in $u>0$, where $(T_u)_{u\ge }$ is the Ornstein-Uhlenbeck semi-group. Here
$\xi$ is not necessarily Malliavin differentiable in the $\mD_{1,1}$-sense.
This provides for Theorem 4.22 of the aforementioned reference an extension to the borderline case $p=1$.
\end{abstract}

\maketitle



\section{Introduction}

Malliavin differentiability and fractional Malliavin differentiability on Gaussian Hilbert spaces is  well-studied
in stochastic analysis and has various applications. There are several ways to introduce Malliavin differentiability: For example,
in the  $L_2$-context one can use It\^o's chaos expansion, one can use smooth random variables and closable operators, and one can
use the Ornstein-Uhlenbeck semi-group. A further method, introduced in \cite{Geiss:Ylinen:2021},  is based on a coupling technique.
This coupling technique was originally designed to investigate the regularity of Backward Stochastic Differential Equations (BSDEs)
(see \cite{GGG:12,CGeiss:Steinicke:16}) and
is meanwhile extended to path-dependent SDEs in \cite{Geiss:Zhou:25}. The coupling method from \cite{Geiss:Ylinen:2021} allows for the introduction of a general class
of Besov spaces $\mB_p^\Phi$, including the Malliavin Besov spaces
obtained by real interpolation, but
also new ones, called $\mB_p^{\Phi_\alpha}$. The spaces  $\mB_p^{\Phi_\alpha}$
differ from the spaces generated by real interpolation and are designed for the purpose
to obtain upper bounds for the $L_p$-path regularity of solutions to BSDEs
in a fully non-Markovian setting (see  \cite{Geiss:Ylinen:2021}).

\bigskip

In this article we concentrate on the spaces $\mB_1^{\Phi_2}$. Let us fix a finite time horizon $T>0$ and assume an $m$-dimensional Brownian motion $w=(w(t))_{t\in [0,T]}$,
$w(t):\Omega \to \mR^m$, and a copy $w'(t):\Omega' \to \mR^m$, defined on complete probability spaces
$(\Omega,\cF,\mP)$ and $(\Omega',\cF',\mP')$, respectively, where $\cF=\sigma(w(t): t\in [0,T])\vee \cN$
and  $\cF'=\sigma(w'(t): t\in [0,T])\vee \cN'$ with $\cN$ and $\cN'$ being the $\mP$ and $\mP'$ null sets, respectively, and where
both Brownian motions start identically in zero and have only continuous paths.
For $0\le a < b \le T$ we introduce the coupled Brownian motion $w^{(a,b]}$ on
$(\overline{\Omega},\overline{\cF},\overline{\mP})$, the completed
product of $(\Omega,\cF,\mP)$ and $(\Omega',\cF',\mP')$, by
$$ w^{(a,b]}(t)
   := \begin{cases}
      w(t)                             &: t\in [0,a], \\
      w(a) + w'(t) - w'(a)             &: t\in [a,b], \\
      w(a) + w'(b) - w'(a) + w(t)-w(b) &: t\in [b,T].
      \end{cases}$$
Within $(\overline{\Omega},\overline{\cF},\overline{\mP})$ we let $(\cF^0_t)_{t\in [0,T]}$, $(\cF_t^{(a,b]})_{t\in [0,T]}$, and
$(\overline{\cF}_t)_{t\in [0,T]}$,
be the augmented natural filtrations of
$w$, $w^{(a,b]}$, and $\overline{w}=(w,w')$, respectively.
\medskip

Following
\cite[Chapter 3]{Geiss:Ylinen:2021}, we obtain a bijective isometry
(with respect to the metric \eqref{eqn:metric_L0})
$$ C^{(a,b]}: L_0(\overline{\Omega},\cF^0_T,\overline{\mP})
\to L_0(\overline{\Omega},\cF_T^{(a,b]},\overline{\mP}). $$
Moreover, we extend any $\xi\in L_0(\Omega,\cF,\mP)$ canonically to
$(\overline{\Omega},\overline{\cF},\overline{\mP})$ and denote this extension by $\xi$ as well (there will be no risk of confusion)
and use the notation
$$ \xi^{(a,b]}:=C^{(a,b]} \xi.$$
The following spaces have been introduced in \cite[Definition 4.12, Section 4.6]{Geiss:Ylinen:2021}:

\bd
For $p\in [1,\infty)$, $\alpha>0$, and $\xi \in L_p(\Omega,\cF,\mP)$  we let
$\xi \in \mB_p^{\Phi_\alpha}$ provided that
$$ |\xi|_{\Phi_\alpha,p}:= \sup_{0 \le a < b \le T}\frac{\| \xi-\xi^{(a,b]}\|_{L_p}}{(b-a)^\frac{1}{\alpha}}<\infty
   \sptext{1}{and set}{1}
   \|\xi\|_{\mB_p^{\Phi_\alpha}}:= \left ( \|\xi\|_{L_p}^p + |\xi|_{\Phi_\alpha,p}^p \right )^\frac{1}{p}.$$
\ed
\medskip

These spaces, although originally introduced to estimate the $L_p$-path regularity of BSDEs, give an insight into structural
properties of the Malliavin derivative. In \cite[Theorem 4.22]{Geiss:Ylinen:2021} the following is shown:
If $p\in [2,\infty)$ and $\xi \in \mD_{1,2}\cap L_p$, then for the Malliavin derivative $D\xi$ it holds
that
$$  |\xi|_{\Phi_2,p} \sim \sup_{0\< a < b\< T} \left \| \left ( \frac{1}{b-a} \int_a^b |D_s \xi|^2  ds\right)^\frac{1}{2} \right \|_{L_p}, $$
whereas for $p\in (1,2)$ one has
\begin{equation}\label{eqn:D_vs_coupling_(1,2)}
           \mathop{\rm esssup}_{s\in [0,T]} \left \| D_s \xi\right \|_{L_p}
  \lesssim |\xi|_{\Phi_2,p}
  \lesssim \sup_{0\< a <  b\< T} \left \| \left ( \frac{1}{b-a} \int_a^b |D_s \xi|^2 ds \right )^\frac{1}{2} \right \|_{L_p}.
\end{equation}
The case $p=1$ was not covered by \cite[Theorem 4.22]{Geiss:Ylinen:2021}. 
The main result of this note is an extension to the case $p=1$:
\bigskip

\bt\label{statement:intro:D_vs_B1_2_new}
For $\xi \in \mB_1^{\Phi_2} \cap L_2$  one has
\begin{equation}\label{eqn:statement:intro:D_vs_B1_2_new}
 \sup_{u>0}\, \mathop\esssup_{t\in [0,T]} \,\,\| \bD_t(T_u \xi) \|_{L_{1,\infty}(\ln L)^{-1}} \le  c_{\eqref{eqn:statement:intro:D_vs_B1_2_new}} | \xi |_{\Phi_2,1},
\end{equation}
where $\bD:\mD_{1,2}\to L_2([0,T]\times \Omega;\mR^m)$ is the Malliavin derivative,
$(T_u)_{u\ge 0}$ the Ornstein-Uhlenbeck semi-group, 
$L_{1,\infty}(\ln L)^{-1}$ is defined in \eqref{eqn:LZ-space},
and $c_{\eqref{eqn:statement:intro:D_vs_B1_2_new}}>0$  is an absolute
constant.
\et
\medskip

Note that $\esssup$ does not require any measurability of $t\mapsto \| \bD_t(T_u \xi) \|_{L^{1,\infty}(\ln L)^{-1}}$.
\medskip

\bx\rm
Let $T=m=1$ and $K\in \mR$. From \cite[Example 2.2]{Geiss:Toivola:2009} we know that
$$ 1_{[K,\infty)}(w(1)) \in (L_2(\mR,\gamma),D_{1,2}(\mR\gamma))_{\frac{1}{2},\infty} $$
where $\gamma$ is the standard Gaussian measure on $\mR$ and
$(L_2(\mR,\gamma),D_{1,2}(\mR\gamma))_{\frac{1}{2},\infty}$ denotes
the real interpolation space.
According to \cite[Theorem 4.16]{Geiss:Ylinen:2021} this is  equivalent to 
$$ \| 1_{[K,\infty)}(w(1)) -  1_{[K,\infty)}(w^r(1)) \|_{L_2} \le c \sqrt{r}
  \sptext{1}{for}{1} r\in (0,1) $$
with $w^r(1):= \sqrt{1-r^2}w(1) + r w'(1)$.
As
$$ {\rm cov}(w(1),w^{(a,b]}(1)) = 1-(b-a)
   \sptext{1}{and}{1}
   {\rm cov}(w(1),w^r(1)) = \sqrt{1-r^2}, $$
and
$\mE | 1_{[K,\infty)}(w(1)) -  1_{[K,\infty)}(w^r(1)) \|_{L_2}$ and
$\mE | 1_{[K,\infty)}(w(1)) -  1_{[K,\infty)}(w^r(1)) \|_{L_2}$
depend on these covariances only, we conclude (by a simple computation) that 
 $$ \| 1_{[K,\infty)}(w(1)) -  1_{[K,\infty)}(w^{(a,b]}(1)) \|_{L_2} \le d \sqrt[4]{b-a} \sptext{1}{for all}{1}
 b-a\in (0,1).$$
Hence 
$$ \mE | 1_{[K,\infty)}(w(1)) -  1_{[K,\infty)}(w^{(a,b]}(1)) | \le d^2 \sqrt{b-a}
\sptext{1}{and}{1} 1_{\{w(1)>K\}}  \in \mB_1^{\Phi_2}.$$
Note that $1_{[K,\infty)}(w(1))\not \in \mD_{1,1}$, i.e.
$1_{[K,\infty)}(w(1))$ is not Malliavin differentiable, see \cite[Proposition 1.2.6]{Nualart:2006}.
\ex
\bigskip
In
\cref{sec:proof_statement:intro:D_vs_B1_2_new}
we prove \cref{statement:intro:D_vs_B1_2_new} and
close with \cref{sec:application_BSDEs} by an application 
to BSDEs.
\medskip

{\bf Some notation.}
We denote by $\ln$ the logarithm with respect to the natural basis $e$.
For a probability space $(\Om,\cF, \mP)$, $L_0(\Om, \cF, \mP)$ stands for all equivalence classes of real-valued random variables, which
becomes a complete metric space under the metric
\begin{equation}\label{eqn:metric_L0}
d(\xi,\eta):=\mE\left[\frac{|\xi-\eta|}{1+|\xi-\eta|}\right].
\end{equation}
 Given a random variable
$f:\Omega\to \mR$, we use the weak-$L_1$ and the  Lorentz-Zygmund space
\begin{align}
     \| f \|_{L_{1,\infty}}
&:= \sup_{\lambda>0} \lambda \mP(|f|>\lambda),\notag \\
\|\xi\|_{L_{1,\infty}(\ln L)^{-1}}&:= \inf\left \{ c >0 : \mP(|\xi|> c \lambda)\le e  \frac{\ln \lambda}{\lambda} :  \lambda >e\right \}.
\label{eqn:LZ-space}
\end{align}


\section{Proof of \cref{statement:intro:D_vs_B1_2_new}}
\label{sec:proof_statement:intro:D_vs_B1_2_new}
For general presentations of Malliavin calculus and related topics we refer to the monographs \cite{Nualart:2006}
and \cite{Janson:2009}. We shall use It\^o's chaos decomposition
$$ L_2 = \bigotimes_{n=0}^\infty \cH_n $$
and the Ornstein-Uhlenbeck semi-group $(T_u)_{u\ge 0}$ with
$$ T_u:L_2(\Omega,\cF,\mP)\to L_2(\Omega,\cF,\mP) \sptext{1}{given by}{1}
   T_u: \sum_{n=0}^\infty I_n(f_n) \to \sum_{n=0}^\infty e^{-n u} I_n(f_n).$$
Then we define the Malliavin Sobolev space
$$ \mD_{1,2}:= \left \{ \xi= \sum_{n=0}^\infty I_n(f_n) \in L_2(\Omega,\cF,\mP):
                  \| \xi\|_{\mD_{1,2}}^2 := 
               \sum_{n=0}^\infty (n+1) \| I_n(f_n) \|_{L_2}^2 < \infty
               \right \}.$$
One point of view is that the Malliavin derivative acts as a map $D \xi: [0,T] \times \Omega\to \mR^m$ where $\mR^m$ is equipped with the
euclidean norm. For example, if $\xi=g(w(t_1),\ldots,w(t_n))$
with $0\< t_1 < \cdots < t_n\< T$ and $g\in C_0^\infty(\mR^m \times \cdots \times \mR^m)$, then
a representative is
$$ (D\xi)(t,\omega) = \sum_{i=1}^m \nabla_i g(w(t_1),\ldots,w(t_n))
1_{[0,t_i]}(t).$$
We have that $T_u:L_2 \to \mD_{1,2} \subseteq L_2$ for $u>0$ and
$\lim_{u\to 0} \| \xi - T_u \xi\|_{L_2}=0$, so that we have an $L_2$-approximations of $\xi\in L_2$ by elements
of $\mD_{1,2}$. The following lemma explains why this approximation is useful
for us:
\medskip

\bl
For $0\le a < b \le T$ one has that
$$ \mE | T_u \xi - (T_u \xi)^{(a,b]}| \le  \mE | \xi - \xi^{(a,b]}|. $$
\el

\begin{proof}
By an approximation argument we can assume that $\xi=g(w(t_1),\ldots,w(t_n))$ for
$0\le t_1 < \cdots < t_n \le T$ and a bounded Borel function $g:\mR^n\to \mR$.
Let $(\tilde{w},\tilde{w}')$ be  a further copy of the pair of Brownian motions $(w,w')$
defined on $(\tilde{\Omega},\tilde{\cF},\tilde{\mP})$ and
$(\tilde{\Omega}',\tilde{\cF}',\tilde{\mP}')$, respectively.
The Brownian motion $\tilde{w}^{(a,b]}$ is obtained by coupling
$\tilde{w}$ and $\tilde{w}'$ as before.
By Mehler's formula we have
$$ T_u \xi = \mE_{\tilde \mP}g (e^{-u} w + \sqrt{1-e^{-2u}} \tilde w). $$
This implies that
\begin{align*}
&    \mE_{\overline{\mP}} | T_u \xi - (T_u \xi)^{(a,b]}| \\
& =  \mE_{\overline{\mP}} \left |  \mE_{\tilde{\mP}} g (e^{-u} w + \sqrt{1-e^{-2u}} \tilde w)  -
             \mE_{\tilde{\mP}\otimes \tilde{\mP}'} g \left (e^{-u} w^{(a,b]} + \sqrt{1-e^{-2u}} \tilde w^{(a,b]}\right ) \right | \\
&\le \mE_{\overline{\mP}} \mE_{\tilde{\mP}\otimes \tilde{\mP}'} \left | g (e^{-u} w + \sqrt{1-e^{-2u}} \tilde w)  -
                             g \left (e^{-u} w^{(a,b]} + \sqrt{1-e^{-2u}} \tilde w^{(a,b]}\right )\right | \\
& = \mE_{\overline{\mP}}  | \xi (w)  - \xi (w^{(a,b]})|. \qedhere
\end{align*}
\end{proof}
\medskip

We recall Burkholder's weak $L_1$-inequality
\cite{Burkholder;1966,DCCox:1982,Osekowski:2011}:
\medskip

\bl\label{statement:weak:L1_BDG}
Let $M=(M_t)_{t\in [0,T]}$ with  $M_t := \int_0^t L(s) dw(s): \Omega\to \mR$
be a continuous martingale where $L=(L(s))_{s\in [0,T]}$ is progressively measurable and  $\mE \int_0^T |L(s)|^2 ds < \infty$. Then,
for $S_2M := \sqrt{\int_0^T |L(t)|^2 dt}$ and some absolute constant $c_\eqref{statement:weak:L1_BDG}>0$, one has that
$$ \| S_2 M\|_{L_{1,\infty}} \le c_\eqref{statement:weak:L1_BDG} \mE |M_T|.$$
\el

\br\rm
There is a direct argument for this inequality:
We use  the good $\lambda$-inequality  for the square function $S_2$ and the maximal function $M^*$,
$$ \mP(S_2M > \beta \lambda, M^* \le \delta \lambda)
   \le \epsilon \mP(S_2M > \lambda)
\sptext{1}{for all}{1} \lambda >0,
$$
for $\beta>1$, $\delta>0$, and $\epsilon:=\frac{\delta^2}{\beta^2-1}>0$,
where $\delta>0$ is chosen such that
$\beta \epsilon <1$, see \cite[Inequality (9.1)]{Burkholder:1973}. From
this we get
\begin{align*}
      \mP(S_2M > \beta \lambda)
&\le  \mP(S_2M > \beta \lambda, M^* \le \delta \lambda)
      +  \mP(M^* > \delta \lambda)\\
&\le  \epsilon \mP(S_2M > \lambda)
      +  \mP(M^* > \delta \lambda)
\end{align*}
and
$
\frac{1}{\beta} \left \| S_2M \right \|_{L_{1,\infty}}
\le \varepsilon \left \| S_2M \right \|_{L_{1,\infty}} +
\frac{1}{\delta} \left \| M^* \right \|_{L_{1,\infty}}
$. Therefore, by Doob's maximal inequality,
$$
\left ( \frac{1}{\beta} - \varepsilon \right ) \left \| S_2M \right \|_{L_{1,\infty}}
\le \frac{1}{\delta} \left \| M^* \right \|_{L_{1,\infty}}
\le \frac{1}{\delta} \left \| M_T \right \|_{L_1}. \qedhere
$$

\er
\medskip

\begin{proof}[Proof of \cref{statement:intro:D_vs_B1_2_new}]
(a) For notational convenience we use for $T_u \xi$ the notation $\xi$.
We start with \cite[Step (b) of the proof of Theorem 4.22]{Geiss:Ylinen:2021}:
For any $0\le a < b \le T$ relation (22) of \cite{Geiss:Ylinen:2021} turns into
$$ \left \| \xi - \xi^{(a,b]}\right \|_{L_1}
   \ge \left \| \xi - \mE [\xi| \cG_a^b] \right \|_{L_1}
   \ge \frac{1}{c_\eqref{statement:weak:L1_BDG}} \left \| \left ( \int_a^b |\mu_s^b|^2 ds \right )^\frac{1}{2} \right \|_{L_{1,\infty}}, $$
where the factor $2$ in the first inequality of (22) is not needed in this direction (see \cite[Lemma 4.23]{Geiss:Ylinen:2021}),
and in the second inequality we use 
\cref{statement:weak:L1_BDG}. More precisely, 
\cref{statement:weak:L1_BDG} is used as replacement for the Burkholder-Davis-Gundy inequality on line 7 of page 103 of \cite{Geiss:Ylinen:2021} which enters (22).
Fixing $0\le a < b \le T$, for $n\ge 0$ we let
$$ t_k^n := a + \frac{k}{2^n}(b-a)
   \sptext{1}{for}{1}
   k=0,\ldots, 2^n $$
and
$$ b_n(s) := \inf \{ t_k^n: s\le t_k^n, k=0,\ldots,2^n \}
   \sptext{1}{for}{1} s\in [a,b].$$
With this we get
\begin{align*}
(b-a) |\xi|^2_{\Phi_2,1}
& = \sum_{k=1}^{2^n} (t_k^n - t_{k-1}^n) |\xi|^2_{\Phi_2,1} \\
&\ge \sum_{k=1}^{2^n} \left \| \xi - \mE [\xi| \cG_{t_{k-1}^n}^{t_k^n}] \right \|_{L_1}^2 \\
&\ge \frac{1}{c^2_\eqref{statement:weak:L1_BDG}}  \sum_{k=1}^{2^n}
     \left \| \left ( \int_{t_{k-1}^n}^{t_k^n} |\mu_s^{b_n(s)}|^2 ds \right )^\frac{1}{2} \right \|^2_{L_{1,\infty}}.
\end{align*}
Hence, for $p\in (0,1)$,
\begin{align}
     \sqrt{b-a} |\xi|_{\Phi_2,1}
&\ge \frac{1}{c_\eqref{statement:weak:L1_BDG}}  \left ( \sum_{k=1}^{2^n}
     \left \| \left ( \int_{t_{k-1}^n}^{t_k^n} |\mu_s^{b_n(s)}|^2 ds \right )^\frac{1}{2} \right \|^2_{L_{1,\infty}}\right )^\frac{1}{2} \notag \\
&\ge  \frac{(1-p)^\frac{1}{p}}{c_\eqref{statement:weak:L1_BDG}}
     \left ( \sum_{k=1}^{2^n}
     \left \| \left ( \int_{t_{k-1}^n}^{t_k^n} |\mu_s^{b_n(s)}|^2 ds \right )^\frac{1}{2} \right \|_{L_p}^2\right )^\frac{1}{2} \label{eqn:needed_inequality} \\
&\ge \frac{(1-p)^\frac{1}{p}}{c_\eqref{statement:weak:L1_BDG}}
     \left ( \sum_{k=1}^{2^n}
     \left [ \left ( \int_{t_{k-1}^n}^{t_k^n} \|\mu_s^{b_n(s)}\|_{L_p}^2 ds \right )^\frac{1}{2}\right ]^2\right )^\frac{1}{2} \notag \\
& = \frac{(1-p)^\frac{1}{p}}{c_\eqref{statement:weak:L1_BDG}}
     \left ( \int_a^b \|\mu_s^{b_n(s)}\|_{L_p}^2 ds \right )^\frac{1}{2}. \notag
\end{align}
So we have that
$$  \left ( \frac{1}{b-a} \int_a^b \|\mu_s^{b_n(s)}\|_{L_p}^2 ds \right )^\frac{1}{2}
    \le  \frac{c_\eqref{statement:weak:L1_BDG}}{(1-p)^\frac{1}{p}} |\xi|_{\Phi_2,1}.$$
As in \cite[p. 40]{Geiss:Ylinen:2021} we find a subsequence $(n_l)_{l=1}^\infty$  such that
$$ \lim_{l\to \infty} \mu_s^{b_{n_l}(s)} = \bD_s \xi \quad \lambda\otimes \mP \mbox{-a.e.}$$
Fatou's lemma implies that
$$  \left ( \frac{1}{b-a} \int_a^b \left \|\liminf_{l\to \infty} \left |\mu_s^{b_{n_l}(s)}\right |\right \|_{L_p}^2 ds \right )^\frac{1}{2}
    \le  \frac{c_\eqref{statement:weak:L1_BDG}}{(1-p)^\frac{1}{p}}|\xi|_{\Phi_2,1}.$$
This implies that
$$  \left ( \frac{1}{b-a} \int_a^b \left \| \bD_s \xi\right \|_{L_p}^2 ds \right )^\frac{1}{2}
    \le  \frac{c_\eqref{statement:weak:L1_BDG}}{(1-p)^\frac{1}{p}} |\xi|_{\Phi_2,1}.$$
(b) Therefore there is a null-set
$N^{(p)} \in \mathcal{B}([0,T])$ such that for all $t\in [0,T]\setminus N^{(p)}$ one has
$$
(1-p)^{1/p} \| \bD_t \xi \|_{L_p}
\le c_\eqref{statement:weak:L1_BDG} | \xi |_{\Phi_2,1}.$$
Let $N:= \bigcup_{p\in (0,1)\cap Q} N^{(p)}$. Hence for all $t\in [0,T]\setminus N$ one has
\begin{equation}\label{eqn:Lp-inequality}
    \sup_{p\in Q\cap (0,1)} (1-p)^{1/p} \| \bD_t \xi \|_{L_p}
\le c_\eqref{statement:weak:L1_BDG} | \xi |_{\Phi_2,1}.
\end{equation}
Let $\lambda>e$ be such that $p:=1-\frac{1}{\ln \lambda}\in Q\cap (0,1)$ and $t\in [0,T]\setminus N$.
Using \eqref{eqn:Lp-inequality} we deduce
$$ \mP\left (\frac{|\bD_t\xi|}{c_\eqref{statement:weak:L1_BDG} | \xi |_{\Phi_2,1}} > \lambda\right ) \le
e \frac{\ln \lambda}{\lambda}.$$
As the RHS is continuous in $\lambda$ the inequality holds true for all $\lambda>e$.
\end{proof}

\br\rm
Let $\Omega := [0,1)$ be equipped with the Borel $\sigma$-algebra and the Lebesgue measure.
For $k\ge 1$ we define $A_k:=[1-2^{-k+1},1-2^{-k})$ and set
$$ f_k := 2^k 1_{A_k}.$$
Then we compute
$$ \| f_k \|_{L_{1,\infty}}=1, \quad
    \left ( \sum_{k=1}^n \| f_k \|^2_{L_{1,\infty}}\right )^\frac{1}{2} =\sqrt{n},
   \sptext{1}{and}{1}
   \left (\sum_{k=1}^n |f_k|^2 \right )^\frac{1}{2} = \sum_{k=1}^n 2^k 1_{A_k}.$$
Moreover, we have
$$ \lambda \mP \left ( \left (\sum_{k=1}^n |f_k|^2 \right )^\frac{1}{2} > \lambda \right ) \le 2. $$
Hence
$$  \left ( \sum_{k=1}^n \| f_k \|^2_{L_{1,\infty}}\right )^\frac{1}{2} \ge \frac{\sqrt{n}}{2} \left \| \left (\sum_{k=1}^n |f_k|^2 \right )^\frac{1}{2} \right \|_{L_{1,\infty}}. $$
In other words, we do not have a uniform inequality
$$  \left ( \sum_{k=1}^n \| f_k \|^2_{L_{1,\infty}}\right )^\frac{1}{2} \le c \left \| \left (\sum_{k=1}^n |f_k|^2 \right )^\frac{1}{2} \right \|_{L_{1,\infty}}. $$
This is the reason why we replace in \eqref{eqn:needed_inequality}  the space $L_{1,\infty}$ by $L_p$ for $p\in (0,1)$.
\er


\section{An application to BSDEs}
\label{sec:application_BSDEs}

We let $m=1$ and specialize the setting from \cite[Definition 2.1]{Briand_etal:2003} by considering the BSDE
$$ Y(t) = \xi + \int_t^T f(s,Y(s)) ds - \int_t^T Z(s) d w(s)$$
such that $\xi$ and $f:[0,T]\times \Omega \times \mR\to \mR$
satisfy
\begin{enumerate}
\item $\xi \in L_2$,
\item $ |f(s,\omega,y_1) - f(s,\omega,y_2)|\< L_f |y_1-y_2|$ for all $y_1,y_2\in \mR$ and $\omega \in \Omega$,
\item $(t,\omega)\mapsto f(t,\omega,y)$ is predictable for all $y\in \mR$,
\item $\mE |\int_0^T |f(s,\cdot,0)|ds|^2 < \infty$,
\item $y\mapsto f(s,\omega,y)$ is non-increasing for all $(s,\omega)\in [0,T]\times \Omega$.
\end{enumerate}
We consider the unique $L_2$-solution $((Y(t),Z(t)))_{t\in [0,T]}$, i.e.
$Y$ is continuous and adapted such that $\mE \sup_{t\in [0,T]}|Y(t)|^2 <\infty$ and
$Z$ is predictable such that $\mE\int_0^T |Z(t)|^2 d t< \infty$,
see \cite[Theorem 4.2]{Briand_etal:2003}.
\medskip

\bt
One has 
$$ \mE |Y_t - Y_t^{(a,b]}| \< \mE \left |\xi-\xi^{(a,b]} \right | 
                           + \mE \int_t^T \left | f(s,Y_s^{(a,b]}) - f^{(a,b]}(s,Y_s^{(a,b]})\right | ds.$$
\et
\begin{proof}
We apply the coupling procedure from \cite[Section 6.2, up to line 5 on page 68]{Geiss:Ylinen:2021} and get the two BSDEs
\begin{align*}
Y(t) & = \xi + \int_t^T f(s,Y(s)) ds - \int_t^T Z(s) d w(s), \\
Y^{(a,b]}(t) & = \xi^{(a,b]} + \int_t^T f^{(a,b]}(s,Y^{(a,b]}(s)) ds - \int_t^T Z^{(a,b]}(s) d w^{(a,b]}(s),
\end{align*}
that are seen as BSDEs on the joint stochastic basis 
$(\overline{\Omega},\overline{\cF},\overline{\mP}, (\overline{\cF}_s)_{s\in [0,T]})$.
We take the difference of the equations to get 
\begin{multline*}
 Y(t)- Y^{(a,b]}(t)   = (\xi- \xi^{(a,b]}) + \int_t^T [f(s,Y(s))- f^{(a,b]}(s,Y^{(a,b]}(s))] ds \\
- \left (\int_t^T Z(s) d w(s) - \int_t^T Z^{(a,b]}(s) d w^{(a,b]}(s)\right ). 
\end{multline*}
To this BSDE we apply   \cite[Corollary 2.3]{Briand_etal:2003} for $k=1$, $d=2$, and $p=1$ and get, $\overline{\mP}$-a.s.,
\begin{align*}
| Y(t)- Y^{(a,b]}(t) |
&\< |\xi-\xi^{(a,b]}| + \int_t^T \alpha(s) [f(s,Y(s))- f^{(a,b]}(s,Y^{(a,b]}(s))] ds  \\
& - \left (\int_t^T \alpha(s) Z(s) d w(s) - \int_t^T \alpha(s) Z^{(a,b]}(s) d w^{(a,b]}(s)\right )
\end{align*}
with 
$$\alpha(s):= \frac{Y(s)-Y^{(a,b]}(s)}{|Y(s)-Y^{(a,b]}(s)|} 1_{\{Y(s)\neq Y^{(a,b]}(s)\}}.$$
Writing 
\begin{multline*}  
   f(s,Y(s))- f^{(a,b]}(s,Y^{(a,b]}(s)) \\
=  [f(s,Y(s))- f(s,Y^{(a,b]}(s))] + [f(s,Y^{(a,b]}(s))- f^{(a,b]}(s,Y^{(a,b]}(s))]
\end{multline*}
and using the monotonicity of $f$ in $y$ yields to
$$ |Y_t - Y_t^{(a,b]}| \< \mE\left [ \left |\xi-\xi^{(a,b]} \right | \big | \overline{\cF}_t \right ]
                           + \mE \left [ \int_t^T \left | f(s,Y_s^{(a,b]}) - f^{(a,b]}(s,Y_s^{(a,b]})\right | ds \big | \overline{\cF}_t \right ]. $$
\end{proof}

\bx\rm
Assume a predictable $U:[0,T]\times \Omega\to \mR$  such that for all
$0\< a <b \<T$ one has for a.a. $(s,\overline{\omega})\in [0,T]\times \overline{\Omega}$ that 
$$ |f(s,\omega,y)-f^{(a,b]}(s,\overline{\omega},y)| \< |U(s,\omega)-U^{(a,b]}(s,\overline{\omega})|
   \sptext{1}{for all}{1} y\in \mR.  $$
If one has  
$$ \mE \int_0^T |U(s) - U^{(a,b]}(s)| ds \< c \sqrt{b-a},$$
then
$$ \mE |Y_t - Y_t^{(a,b]}| \< \mE \left |\xi-\xi^{(a,b]} \right | 
                           + c \sqrt{b-a}.$$

\ex
\bigskip


{\sc Acknowledgment}. The author thanks Yao Chen, Jiagang Ren, and Hua Zhang for helpful discussions.
\medskip

{\sc Disclaimer}. The author did not not use any AI-tools, neither to obtain the results nor for the  finalization of the manuscript.


\end{document}